\documentclass[11pt]{amsart}

\catcode`\@=11

\long\def\@savemarbox#1#2{\global\setbox#1\vtop{\hsize\marginparwidth 
  \@parboxrestore\tiny\raggedright #2}}
\newcommand\lref[1]{\ref{#1}%
\@ifundefined{r@DisplaY #1}{}{ (#1)}}

\newcommand\fakelabel[2]{\@bsphack\if@filesw {\let\thepage\relax
   \newcommand\protect{\noexpand\noexpand\noexpand}%
\xdef\@gtempa{\write\@auxout{\string
      \newlabel{#1}{{#2}{\thepage}}}}}\@gtempa
   \if@nobreak \ifvmode\nobreak\fi\fi\fi\@esphack}

\catcode`\@=12

\def\Empty{}
\newcommand\oplabel[1]{
  \def\OpArg{#1} \ifx \OpArg\Empty {} \else
        \label{#1}
  \fi}
\newtheorem{theoremSt}{Theorem}[section]

\newtheorem{exampleSt}[theoremSt]{Example}
\newtheorem{exerciseSt}[theoremSt]{Exercise}

\newcommand\MakeStEnv[1]{
  \newenvironment{#1}[1]{
  \begin{#1St} \oplabel{##1}%
  \global\def\CrntSt{\thetheoremSt}%
}{ 
  \end{#1St} }
  \newenvironment{#1+}[1]{
  \begin{#1St} \label{##1}%
  \label{DisplaY ##1}%
  \global\def\CrntSt{\thetheoremSt}%
  \def\Labl{##1}\ifx\Labl\Empty{} \else {\em (\Labl)\,}\fi%
}{ 
  \end{#1St} }
}
\MakeStEnv{theorem}
\MakeStEnv{corollary}
\MakeStEnv{proposition}
\MakeStEnv{lemma}
\MakeStEnv{definition}
\MakeStEnv{conjecture}

\long\def\state#1#2{
\medskip\par\noindent
{\bf #1 } 
{\it #2}
\par\medskip
}

\long\def\realfig#1#2#3{
\begin{figure}[htbp]
\centerline{\psfig{figure=#2}}
\caption[#1]{#3}
\oplabel{#1}
\end{figure}}

\newlength{\saveu}

\newenvironment{pf*}[1]{%
 \begin{proof}[#1]%
}{ 
 \end{proof}
}

\newcommand{\finishproof}[1]{ 
  \def\FPArg{#1}
  \ifx\FPArg\Empty
        \newcommand\FPArg{\CrntSt}  \fi
  \smallbreak\noindent\makebox[\textwidth]{\hfill\fbox{\FPArg}}
  \medbreak\noindent
}

\newcommand\AAA{{\mathcal A}}

\newcommand\CC{{\mathcal C}}

\newcommand\FF{{\mathcal F}}
\newcommand\GG{{\mathcal G}}

\newcommand\LL{{\mathcal L}}
\newcommand\MM{{\mathcal M}}

\newcommand\PP{{\mathcal P}}

\newcommand\TT{{\mathcal T}}

\newcommand\PMF{{\PP\kern-2pt\MM\FF}}
\newcommand\ML{{\MM\LL}}
\newcommand\PML{{\PP\kern-2pt\MM\LL}}
\newcommand\GL{{\GG\LL}}

\newcommand\ep{\epsilon}

\newcommand\hhat{\widehat}

\newcommand\union{\cup}
\newcommand\intersect{\cap}
\newcommand\bbR{{\mathord{\text{I\kern-2pt R}}}}        
\newcommand\bbH{{\mathord{\text{I\kern-2pt H}}}}        

\newcommand\C{{\mathbf C}}

\newcommand\Z{{\mathbf Z}}
\newcommand\R{{\mathbf R}}
\newcommand\Q{{\mathbf Q}}

\newcommand\Hyp{{\mathbf H}}

\newcommand\PSL[1]{\text{PSL}_{#1}}

\newcommand\SL[1]{\text{SL}_{#1}}

\newcommand\bigrightarrow[1]{\hbox to #1{\rightarrowfill}}
\newcommand\bigleftarrow[1]{\hbox to #1{\leftarrowfill}}

\newcommand\boundary{\partial}
\newcommand\semidir{\mathrel{\hbox{\vrule depth-.03ex height1.1ex\kern-0.15em$\times$}}}

\newcommand{\diam}{\operatorname{diam}}

\renewcommand{\Im}{\operatorname{Im}}

\numberwithin{equation}{section}

\input{psfig}

\newcommand{\T}{{\mathbf T}}

\newcommand{\pleat}{\operatorname{\mathbf{pleat}}}
\newcommand{\short}{\operatorname{\mathbf{short}}}

\newcommand{\collar}{\operatorname{\mathbf{collar}}}

\newtheorem{question}[theoremSt]{Question}

\begin{document}

\title{Short geodesics and end invariants}
\author{Yair N. Minsky}
\address{SUNY Stony Brook}
\date{March 27, 2000}

\maketitle

\renewcommand\marginpar[1]{} 

Even topologically simple hyperbolic 3-manifolds
can have very intricate geometry. Consider in particular a closed
surface $S$ of genus 2 or more, and the product $N=S\times\R$. This
3-manifold admits a large family of complete, infinite-volume
hyperbolic metrics, corresponding to faithful representations
$\rho:\pi_1(S)\to \PSL 2(\C)$ with discrete image. 

The geometries of $N$ are very different from the product structure that
its topology would suggest. Typically, $N$ contains a complicated
pattern of ``thin'' and ``thick'' parts. The thin parts are collar
neighborhoods of very short geodesics, typically infinitely many. Each
one, called a ``Margulis tube'', has a well-understood shape, but the
way in which these are arranged in $N$, and in particular the
identities of the short geodesics as elements of the fundamental
group, are still something of a mystery.

This issue is closely related to the basic classification conjecture
associated with these manifolds, Thurston's ``ending lamination
conjecture''. This  conjecture states that certain asymptotic
invariants of the 
geometry of $N$, called ending invariants, in fact determine $N$ completely.
(Actually the classification of hyperbolic structures for any manifold
with incompressible boundary reduces to this case, by restriction to
boundary subgroups.)

In this expository paper we will focus on the following question: What
information do the ending invariants give about the presence of very
short geodesics in the manifold? We will summarize and discuss 
the theorem below, part of whose proof appears in \cite{minsky:kgcc}
and part of which will be in \cite{minsky:bdgeom}, as well as a few
conjectures. 

\state{Bounded Geometry Theorem.}{Let $S$ be a closed surface, and 
consider a Kleinian surface group 
$\rho:\pi_1(S)\to \PSL 2(\C)$ with
no externally short curves, and ending invariants $\nu_+$ and $\nu_-$. 
Then
$$
\inf_{\gamma\in\pi_1(S)} \ell_\rho(\gamma) > 0 
\iff
\sup_{Y\subset S} d_Y(\nu_+,\nu_-) < \infty.
$$
}

Here the supremum is over proper essential isotopy classes of
subsurfaces in $S$, and
the quantities $d_Y(\nu_+,\nu_-)$, called ``projection
coefficients'', are defined in Section \ref{define dY}. The quantity
$\ell_\rho(\gamma)$ is the translation distance of $\rho(\gamma)$ in
$\Hyp^3$, or
the length of the closed geodesic associated to $\gamma$ in the 3-manifold.
(The condition on externally short curves is not really necessary -- 
it is added to simplify the other definitions and discussions -- see
\S\ref{surface groups} below).

Part of our goal is to advertise a combinatorial object known as the
{\em complex of curves on a surface},  as a tool for studying the
geometry of hyperbolic 3-manifolds. This object 
is used for definining the coefficients $d_Y$, and 
in general it encodes something about 
the structure of the set of simple loops on a surface. 
In particular, face transitions between simplices in this complex
correspond to elementary moves on pants decompositions of $S$, and
these in turn correspond to homotopies between elementary pleated
surfaces in a hyperbolic 3-manifold. The interaction between the
combinatorial and geometric aspects of these moves is our main object
of study, and seems to be worthy of further consideration.

\section{Definitions}
\label{defs}
\subsection{Surface groups and  ending laminations}
\label{surface groups}

Let $S$ be a closed surface of genus $g\ge 2$.
A {\em Kleinian surface group} will be a
representation $\rho:\pi_1(S)\to \PSL 2(\C)$, discrete and faithful.
The quotient $\Hyp^3/\rho(\pi_1(S))$ is denoted $N_\rho$, and comes
equipped with a homotopy class of homotopy equivalences $S\to N_\rho$,
determined by $\rho$. In fact $N_\rho$ is homeomorphic to $S\times\R$,
by Thurston's theory of tame ends \cite{wpt:notes} and Bonahon's
Tameness theorem \cite{bonahon}. 

We can associate to $\rho$ two {\em ending invariants} $\nu_-$ and $\nu_+$,
which we will describe in the special case that $\rho$ has {\em no
  parabolics} (see also \cite{minsky:knoxville} and Ohshika
\cite{ohshika:ending-lams}).

Let $C(N_\rho)$ be the {\em convex core} of $N_\rho$, the smallest convex
submanifold whose inclusion is a homotopy equivalence. Fixing an
orientation on $S$ and $N_\rho$, there is an orientation-preserving
homeomorphism of $N_\rho$ to $S\times\R$ taking
$C(N_\rho)$ onto exactly one of 
$S\times\R$, $S\times [0,\infty)$, $S\times(-\infty,1]$ or
$S\times[0,1]$. 

The end of $N$ defined by neighborhoods $S\times(a,\infty)$ is called
$e_+$, and the one defined by $S\times(-\infty,a)$ is called
$e_-$. If an end's neighborhoods all meet the convex hull it
is called {\em geometrically infinite}, and otherwise it is {\em
  geometrically finite}. Suppose $e_+$ is geometrically
finite. Then the component $\boundary_+(C(N_\rho))$ corresponding to
$S\times\{1\}$ is a convex surface, and its exterior
$S\times(1,\infty)$ develops out to a ``conformal structure at
infinity'' on $S$, which we call $\nu_+$. (This surface is obtained
from the action of $\rho(\pi_1(S)$ on the Riemann sphere). We define
$\nu_-$ in the same way when $e_-$ is geometrically finite.

Thurston pointed out that boundary $\boundary_+(C(N_\rho))$ is itself
a hyperbolic surface; let us call its structure $\nu'_+$. A theorem of
Sullivan (proof in Epstein-Marden \cite{epstein-marden}) states that
$\nu'_+$ and $\nu_+$ differ by a  uniformly bilipschitz distortion.

To describe the invariant of a geometrically infinite end we need to
briefly recall the notion of a {\em geodesic lamination}. Fixing a
hyperbolic metric on $S$, a
geodesic lamination is a closed subset of $S$ foliated by geodesics. 
Let $\GL(S)$ denote the set of all of these. A {\em measured
lamination} is a geodesic lamination equipped with a Borel measure on
transverse arcs, invariant under transverse isotopy. The space $\ML(S)$
of measured laminations admits a natural topology coming from  weak-*
convergence of the measures. On the supporting geodesic laminations, 
this is related to but not quite the same as the topology of Hausdorff
convergence. However the difference will not be important to us here. 
Simple closed geodesics with positive weights are
dense in $\ML(S)$, and we will consider geodesic laminations obtained
as supports of limits in $\ML(S)$ of sequences of simple closed
curves. Finally we remark that the choice of metric on $S$ is
irrelevant, as any other choice yields naturally isomorphic spaces of
laminations. For more details on this topic see 
Bonahon \cite{bonahon:curves,bonahon:laminations},
Canary-Epstein-Green \cite{ceg}, or Casson-Bleiler
\cite{casson-bleiler}.

If $e_+$ is geometrically infinite then the convex hull contains an
infinite sequence of closed geodesics $\gamma_n$, all homotopic to
{\em simple} closed loops on $S$, and eventually contained in
$S\times(a,\infty)$ for any $a$. This is a theorem of Bonahon, and
Thurston (previously) showed that for such a sequence the curves on
$S$ must converge in the sense of the previous paragraph to a
unique geodesic lamination on $S$. We call this lamination
$\nu_+$, the ending lamination of $e_+$. The
corresponding discussion for $e_-$ gives $\nu_-$.

Finally let us define the technical simplifying condition in the
statement of the Bounded Geometry Theorem.
Call a curve $\gamma$ in $S$ {\em externally short}, with respect to a
representation $\rho$, if it is either parabolic or has length less
than $\ep_1$ with respect to the structures $\nu_-$ and $\nu_+$ (if
these are not laminations),
where $\ep_1$ 
is some fixed constant small enough that there exist hyperbolic
structures on $S$ with no curves of length less than $\ep_1$.
Note in particular that if $\rho$ has two degenerate ends then it
automatically has no externally short curves.

\subsection{Pleated surfaces}

A {\em pleated surface} is a map $f:S\to N$ together with a hyperbolic
 metric on $S$, written $\sigma_f$ and called the {\em induced
 metric}, and a $\sigma_f$-geodesic lamination $\lambda$ on $S$,
 called the pleating locus, so that the following holds: 
 $f$ is length-preserving on paths, maps leaves of $\lambda$ to geodesics, and
 is totally geodesic on the complement  of $\lambda$. These were
 introduced by Thurston \cite{wpt:notes}, and we will see some
 explicit examples in \S \ref{pleat on pants}.

It is a consequence of the work of Thurston and Bonahon that a
geometrically infinite end of a surface group $\rho$ admits pleated
surfaces in the homotopy class of $\rho$ contained in any neighborhood
of the end. The pleating loci of these  surfaces must converge to the
ending lamination, and their hyperbolic structures converge to this
lamination in Thurston's compactification of the Teichm\"uller space. 

\subsection{Complexes of arcs and curves:}
\label{complex defs}

Let $Z$ be a compact finite genus surface, possibly with boundary.
If $Z$ is not an annulus, 
define  $\AAA_0(Z)$ to be 
the set of  essential homotopy classes of 
simple closed curves or properly embedded arcs in $Z$. Here
``homotopy class'' means free homotopy for closed curves, and
homotopy rel $\boundary Z$ for arcs.
``Essential'' means the homotopy class does not contain the constant
map or a map into the boundary. If $Z$ is an annulus, we make the same
definition except that homotopy for arcs is rel endpoints.

We can extend $\AAA_0$ to a simplicial complex
$\AAA(Z)$ by letting a $k$-simplex be any $(k+1)$-tuple
$[v_0,\ldots,v_k]$ with $v_i\in\AAA_0(Z)$ distinct and having
pairwise disjoint representatives. 

Let $\AAA_i(Z)$ denote the $i$-skeleton of $\AAA(Z)$, and let
$\CC(Z)$ denote the subcomplex spanned by vertices corresponding to
simple closed curves. This is the ``complex of curves of $Z$''.

If we put a path metric on $\AAA(Z)$ making every simplex regular
Euclidean of sidelength 1, then it is clearly quasi-isometric to its
$1$-skeleton. It is also quasi-isometric to $\CC(S)$ except in a few
simple cases when $\CC(S)$ has no edges. When $\boundary Z=\emptyset$,
of course $\AAA=\CC$.

It is a nice exercise to compute $\AAA(Z)$ exactly for $Z$ a one-holed
torus, and we leave this to the reader. The answer is closely related
to the Farey graph in the plane -- see \cite{minsky:taniguchi}.

Fix our closed surface $S$ and let 
$\GL(S)$ denote the set of geodesic laminations on $S$ (note that
$\AAA_0(S)=\CC_0(S)$ can identified with a subset of $\GL(S)$).
Let $Y\subset S$ be  a proper essential closed subsurface (all boundary curves
homotopically nontrivial). We have a ``projection map'' 
$$
\pi_Y : \GL(S) \to \AAA(\hhat Y)\union\{\emptyset\}
$$
defined as follows: there is a unique
cover of $S$ corresponding to the inclusion
$\pi_1(Y)\subset\pi_1(S)$, which can be
naturally compactified using the circle at infinity
of the universal cover of $S$ to yield a surface $\hhat Y$ homeomorphic
to $Y$ (remove the limit set of $\pi_1(Y)$ and
take the quotient of the rest). Any lamination $\lambda\in\GL(S)$
lifts to this cover as a collection of closed curves or arcs that have
well-defined endpoints in $\boundary \hhat Y$.
Removing the trivial components, we have a simplex of $\AAA(\hhat Y)$
and we can take, say, its barycenter
(we can also get the empty set if there are no essential components).
A version of this projection also appears in Ivanov \cite{ivanov:subgroupsbook,ivanov:rank}.

If $\beta,\gamma\in\GL(S)$ (in particular in $\CC(S)$ have non-trivial
intersection with $Y$, we denote their ``$Y$-distance'' by:

$$
d_Y(\beta,\gamma) \equiv d_{\AAA(\hhat Y)}(\pi_Y(\beta),\pi_Y(\gamma)).
$$
Note that $\AAA(\hhat Y)$ can be naturally identified with $\AAA(Y)$,
except when $Y$ is an annulus, in which case the pointwise
correspondence of the boundaries matters. 
In the annulus case $d_Y$ measures relative twisting of arcs
determined rel endpoints, and in all other cases we ignore twisting on
the boundary of $\hhat Y$. If $\alpha$ is the core curve of an annulus
$Y$ we will also write 
$$
d_\alpha = d_Y.
$$
See \cite{farb-lubotzky-minsky} for an application of this
construction in the annulus case. 

The complex of curves $\CC(S)$ was first introduced by Harvey
\cite{harvey:boundary}. 
It was applied by Harer \cite{harer:stability,harer:cohomdim} and
Ivanov  \cite{ivanov:complexes1,ivanov:complexes2,ivanov:complexes3}
to study the 
mapping class group of $S$. Similar complexes were  introduced by
Hatcher-Thurston \cite{hatcher-thurston}. Masur-Minsky
\cite{masur-minsky:complex1} proved that $\CC(S)$ is
$\delta$-hyperbolic in the sense of Gromov, and then applied this in
\cite{masur-minsky:complex2} to prove the structural theorems on pants
decompositions that we will use in  Section \ref{proof main}.

\subsection{Projection coefficients}
\label{define dY}

Let us now see how to define the coefficients
$$
d_Y(\nu_+,\nu_-)
$$
which appear in the main theorem, where $\nu_\pm$ are ending
invariants for a surface group.
Using $\pi_Y$ as above, we can already define this 
whenever $\nu_\pm$ are laminations.
In the case of a geometrically finite end when $\nu_+$ or $\nu_-$ are
hyperbolic metrics, we can extend this definition as follows:

If $\sigma $ is a hyperbolic metric on $S$, and $L_1$ a fixed
constant, define
$$
\short(\sigma)$$ 
to be the set of pants decompositions of $S$ with
total $\sigma$-length at most $L_1$. A theorem of Bers 
(see \cite{bers:degenerating,bers:inequality} and
Buser \cite{buser:surfaces})
says that $L_1$ can be chosen,
depending only on genus of $S$, so that $\short(\sigma)$ is always non-empty.
Let us also choose $L_1$ sufficiently large that, if $\sigma$ has no
curves of length less than $\ep_1$ (the constant from the end of
\S\ref{surface groups}), then {\em every} curve in $S$
intersects some $P\in\short(\sigma)$. 

Thus e.g. if both $\nu_+$ and $\nu_-$ are hyperbolic structures, we
may consider distances
$$
d_Y(P_+,P_-)
$$
for any $P_\pm\in\short(\nu_\pm)$ that both intersect $Y$ essentially,
and notice that the numbers obtained cannot vary by more than a
uniformly bounded constant. We let $d_Y(\nu_+,\nu_-)$ be, say,
the minimum over all choices. The case when one of
$\nu_\pm$ is a lamination and the other is a hyperbolic metric is handled
similarly. Note that the condition that $\rho$ has no externally short
curves implies that $d_Y(\nu_+,\nu_-)$ is well-defined for all $Y$.

\section{Margulis tubes}

Let $\gamma$ be a loxodromic element of a Kleinian group $\Gamma$. 
We denote its complex translation length by $\lambda(\gamma) =
\ell+i\theta$ (determined mod $2\pi i$).  Let $\TT_\ep$ be the
$\gamma$-invariant set 
$\{x\in\Hyp^3: \inf_n d(x,\gamma^n(x)) \le \ep\}$. If $\ell(\gamma)<\ep$
This is a tube of some radius $r$ around the axis of $\gamma$, and 
The Margulis Lemma and Thick-Thin decomposition tell us
(see e.g. \cite{kazhdan-margulis,wpt:textbook,benedetti-petronio}) that 
there is a universal constant $\ep_0$ such that if $\ell(\gamma)<\ep_0$
then $\TT_{\ep_0}/\langle\gamma\rangle$ embeds as a solid torus $\T_\gamma$
in  $N=\Hyp^3/\Gamma$, called a {\em Margulis tube}, and furthermore
that all Margulis tubes in $N$ are disjoint.

The radius $r$ of the tube goes to $\infty$ as the length of the core
goes to 0. See Brooks-Matelski \cite{brooks-matelski} and Meyerhoff
\cite{meyerhoff:volumes}
for more precise bounds.

Thus in some sense the geometry around a very short curve in $N$ is
very well understood. It is more difficult to determine the
pattern in which these tubes are arranged in the manifold, and in
particular which curves $\gamma$ have length less than a given $\ep$.

\subsection{Margulis tubes in surface groups}
\label{margulis tubes}

When $\Gamma$ is the image $\rho(\pi_1(S))$ of a Kleinian surface
group, there is a little more we can say. An  observation of Thurston
\cite{wpt:II}, together with Bonahon's tameness theorem
\cite{bonahon}, imply that {\em only simple curves can be short:} that is,
$\ep_0$ may be chosen so that,
if $\ell_\rho(\gamma)<\ep_0$
and $\gamma$ is a primitive element of $\pi_1(S)$ then $\gamma$ is
represented by a simple loop in $S$. This is because, by Bonahon's
theorem, every point in $N_\rho$ is uniformly near the image of a
pleated surface. Thurston pointed out using a simple area bound that
if $\ep_0$ is sufficiently short a
$\pi_1$-injective pleated surface can only meet $\T_\gamma$
in the image of its own 2-dimensional Margulis
tube. The core of this tube must therefore be $\gamma$. 


\subsection{Bounds}

An {\em upper bound} on the length of a curve in a surface group can
be obtained in terms of the conformal boundary at infinity. Bers
showed \cite{bers:boundaries} for a Quasi-Fuchsian representation $\rho$, that
$$
\frac1{\ell_\rho(\gamma)} \ge \frac12\left(\frac1{\ell_+(\gamma)} +
  \frac1{\ell_-(\gamma)}\right) 
$$
where $\ell_\pm$ denote lengths in the hyperbolic
structures on $S$ coming from the two conformal structures $\nu_\pm$
at infinity. The argument uses a monotonicity property for conformal
moduli and the action of $\gamma$ on the Riemann sphere. 
When $S$ is a once-punctured  torus this upper bound can be
generalized to an estimate in both directions (see
\cite{minsky:torus}). In general we have no such result, but in
Section \ref{conjectures} we will state a conjectural estimate.


\section{Bounded geometry}

We say that $\rho$ has {\em bounded geometry} if there is a positive lower
bound on the translation lengths of all group elements. 
This condition incidentally disallows parabolic elements (in a more
general discussion we would allow them and revise the condition),
but the real point is 
that there is a positive lower bound on the lengths of all closed
geodesics in the quotient manifold.

In \cite{minsky:slowmaps,minsky:endinglam}, we showed that 
bounded geometry implies a positive solution
to the ending lamination conjecture. That is, if $\rho_1$ and $\rho_2$
both have bounded geometry, and have the same ending invariants, then
they are conjugate in $\PSL 2(\C)$. This result was accompanied by a
fairly explicit bilipschitz model for the metric on $N$, derived from
the Teichm\"uller geodesic joining the two ending laminations.

The Bounded Geometry theorem gives us a way to strengthen
this result,  since it implies that bounded geometry is detected by
the ending invariants:

\begin{corollary}{Bounded Geometry ELC}
Let $\rho_1,\rho_2$ be Kleinian surface groups with the same ending
invariants, and suppose that $\rho_1$ has bounded geometry. Then
$\rho_1$ and $\rho_2$ are conjugate in $\PSL 2(\C)$.
\end{corollary}

It is worth noting that bounded geometry is a rare condition.  In the
boundary of a Bers slice, for example, there is a topologically
generic (dense $G_\delta$) set of representations each of which has
arbitrarily short elements (see McMullen
\cite[Cor. 1.6]{mcmullen:cusps}, and Canary-Culler-Hersonsky-Shalen
\cite{canary-culler-hersonsky-shalen:cusps} for generalizations).

%
%
%
%
%
%
%

\section{The proof of the bounded geometry theorem}
\label{proof main}


The proof of the direction $(\Rightarrow)$ of 
the Bounded Geometry Theorem appears in
\cite{minsky:kgcc}. The essential tool used there is Thurston's
``efficiency of pleated surfaces'' theorem from \cite{wpt:II}.
We will outline the proof of $(\Leftarrow)$, for
which the details will appear in \cite{minsky:bdgeom}.

\medskip

In roughest form, the argument is this: Let $\gamma\in\pi_1(S)$ be an
element with $\ell_\rho(\gamma) < \ep_0$, and let $\T_\gamma$ be its
Margulis tube. We will use the condition $\sup d_Y(\nu_+,\nu_-) < \infty$ to
construct a sequence of pleated surfaces $\{f_i\}_{i=0}^M$ with the
following properties:
\begin{enumerate}
\item The size $M$ of the sequence is bounded by 
$$M\le K (\sup d_Y(\nu_+,\nu_-))^a$$
where $K,a$ depend only on the genus of $S$.
\item Any successive $f_i,f_{i+1}$ are connected by a homotopy 
$H:S\times[i,i+1]\to N_\rho$ which is  {\em uniformly bounded}
except in a special case, described below.
\item The total homotopy $H:S\times[0,M]\to N_\rho$ homologically
encloses $\T_\gamma$.
\end{enumerate}

Part (3) means that the image of $H$ must cover all of $\T_\gamma$.
Thus, if the ``special case'' of (2) does not occur, then the bounds of (1) and
(2) give a uniform diameter bound on $\T_\gamma$, and hence a lower
bound on $\ell_\rho(\gamma)$. 

The ``special case'' of (2) corresponds to the curve $\gamma$ itself
appearing in the pleating locus of some subsequence of the $f_i$. In
this case  a more delicate argument is needed, using the 
annulus projection
distance $d_\gamma(\nu_+,\nu_-)$ to bound the size of $\T_\gamma$.

Let us now introduce the ingredients needed for this construction. In
\S\ref{finishing argument} we will return to the main proof.

\subsection{Adapted pleated surfaces}
\label{pleat on pants}

If $Q$ is a collection of disjoint, homotopically distinct curves on
$S$ (henceforth a ``curve system''), and $\rho$ a fixed Kleinian surface group,
we let $\pleat(Q,\rho)$ denote the set of pleated surfaces $f:S\to
N_\rho$, in the homotopy class determined by $\rho$, which map
representatives of  $Q$ to geodesics. 
There is the usual equivalence relation on this set, 
in which $f\sim f\circ h$ if  $h$ is a
homeomorphism of $S$ homotopic to
the identity. Let $\sigma_f$ denote the
hyperbolic metric on $S$ induced by $f$.

In particular, if $Q$ is a maximal curve system, or ``pants
decomposition'', $\pleat(Q,\rho)$ consists of finitely many
equivalence classes, all
constructed as follows: Extend $Q$ to a triangulation of $S$ with
one vertex  on each component of $Q$,
and ``spin'' this triangulation around $Q$,
arriving at a lamination $\lambda$ whose closed leaves are $Q$ and
whose other leaves spiral onto $Q$, as in Figure \ref{spun
lamination}. 

\realfig{spun lamination}{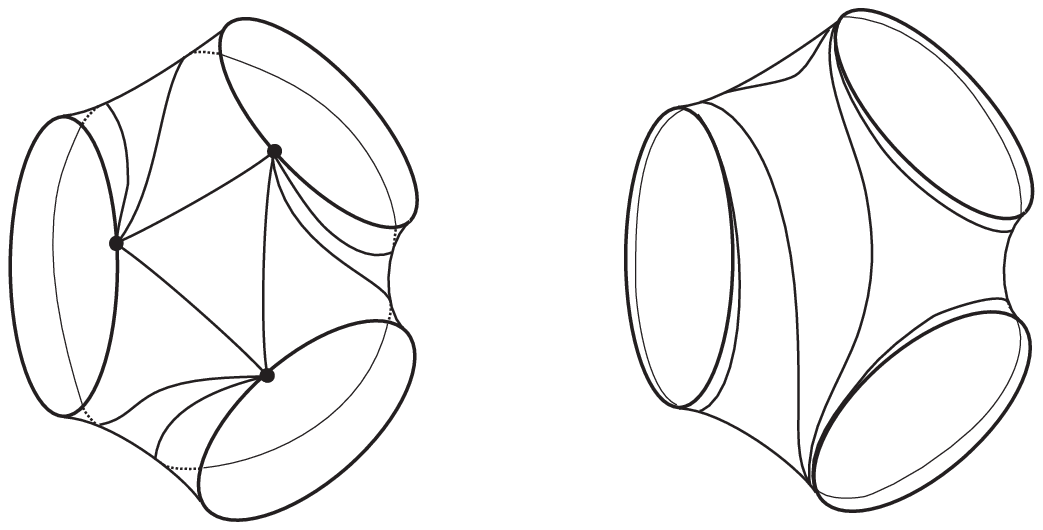}{The lamination obtained by spinning a
triangulation around a curve system. The picture shows one pair of
pants in a decomposition.}

A unique pleated surface (up to equivalence) exists carrying
$\lambda$ to geodesics, since no element of $Q$ is parabolic (by
hypothesis on $\rho$). This was originally observed by Thurston (see
\cite{wpt:notes} and Canary-Epstein-Green \cite[Thm 5.3.6]{ceg} for a
proof).
The choices of $\lambda$ coming from the finite number of possible
triangulations up to isotopy, and the
different directions of spiraling, account for all of $\pleat(Q,\rho)$.

\subsection{Elementary moves}

An elementary move on a maximal curve system $P$ is a replacement of a
component $\alpha$ of $P$ by $\alpha'$, disjoint from the rest of $P$,
so that $\alpha$ and $\alpha'$ 
are in one of the two configurations shown in Figure 
\ref{elementary move fig}.

\realfig{elementary move fig}{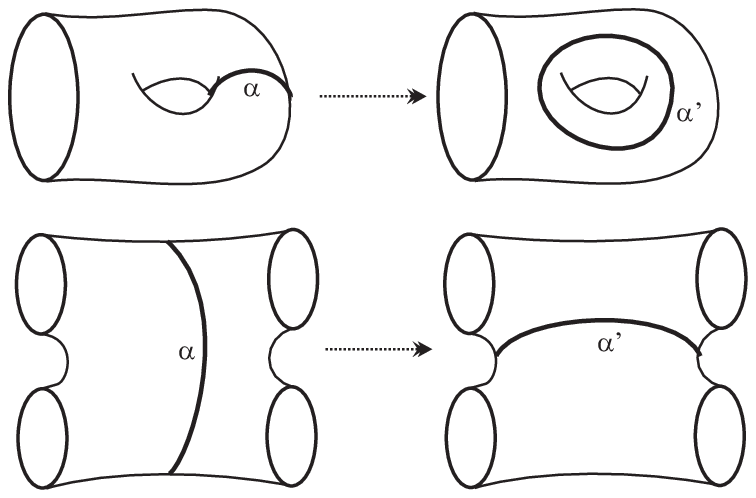}{The two types of elementary moves.}

We indicate this by  $P\to P'$ where
$P'=P\setminus\{\alpha\}\union\{\alpha'\}$ is the new curve system.
Note that there are infinitely many choices for $\alpha'$, naturally
indexed by $\Z$.

Pleated surfaces associated to an elementary move are homotopic in a
controlled way. Let us first recall (see Buser \cite{buser:surfaces})
that a simple geodesic $\gamma$ in a hyperbolic surface $(S,\sigma)$
always admits a ``standard collar'', which is an annulus of radius
depending only on 
$\ell_\sigma(\gamma)$, such that disjoint geodesics have disjoint
collars, and when $\ell_\sigma(\gamma) < \ep_0$ the collar covers all
but a bounded part of the 
$\ep_0$-Margulis tube. We write this collar as 
$\collar(\gamma,\sigma)$, or $\collar(P,\sigma)$ for the union of
collars over a curve system $P$.

\begin{lemma+}{Elementary Homotopy}
If $P_0\to P_1$ is an elementary move exchanging $\alpha_0$ and
$\alpha_1$, $\rho$ is a Kleinian surface
group, and $f_i \in \pleat(P_i,\rho)$ for $i=0,1$, then
there exists a homotopy $H:S\times[0,1]\to N_\rho$ with the
following properties: 
\begin{enumerate}
\item $H_0\sim f_0$ and $H_1\sim f_1$ under the usual equivalence.
\item If $\sigma_i$ is the induced metric of $H_i$ (for $i=0,1$) 
then $\collar(P_j,\sigma_i) = \collar(P_j,\sigma_{1-i})$, for $j=0,1$.

\item The metrics $\sigma_0$ and $\sigma_1$ are $K$-bilipschitz
  except possibly when $\ell_\rho(\alpha_i) < \ep_0$ for $i=0$ or $1$.
  In that case the metrics are locally $K$-bilipschitz 
  outside $\collar(\alpha_0)\union\collar(\alpha_1)$ (or just one
  collar if only one  curve is short in $N_\rho$).

\item The trajectories $H(p\times[0,1])$ are bounded in length by $K$
except 
possibly when $p\in\collar(\alpha_i)$ and $\ell_\rho(\alpha_i)<\ep_0$,
in which case they are bounded outside of $\T_{\rho(\alpha_i)}$. 

\end{enumerate}
The constant $K$ depends only on the genus of $S$.

\end{lemma+}
(Note that $\collar(\alpha_i)$ in (3) and (4) makes sense without
specifying the metric $\sigma_j$, since the two are equal by (2).)

It is worth pointing out that this theorem applies without any
a-priori bounds on the lengths $\ell_\rho(P_i)$. The proof is an
application of Thurston's Uniform Injectivity theorem for pleated
surfaces, and the closely related
Efficiency of Pleated Surfaces \cite{wpt:II} (see also Canary
\cite{canary:schottky}). 
These theorems control the amount kind of bending that can occur in a
pleated surface, and in particular can be used to compare two pleated
surfaces that share part of their pleating locus.

We also remark that part (2) is just for convenience -- it is easy to
arrange by an appropriate isotopy.

\subsection{Resolution sequences}

In \cite{masur-minsky:complex2}, we show the existence of special
sequences of elementary moves that are controlled in terms of the
geometry of the complex of curves, and particularly the projections
$\pi_Y$. First some terminology: if $P_0\to P_1 \to \cdots \to P_n$ is
an elementary-move sequence and $\beta$ is any simple closed curve,
denote
$$J_\beta = \{ i\in[0,n]: \beta\in P_i\}.$$
(Here $\beta\in P$ means $\beta$ is a component of $P$.)
We also denote $J_{\beta_1,\ldots,\beta_k} = \union J_{\beta_i}$.

Note that if $\beta$ is a curve and $J_\beta$ is an interval $[k,l]$,
then the elementary move $P_{k-1}\to P_k$ exchanges some $\alpha$ for
$\beta$, and $P_l\to P_{l+1}$ exchanges $\beta$ for some
$\alpha'$. Both $\alpha$ and $\alpha'$ intersect $\beta$, and we call
them the {\em predecessor} and {\em successor} of $\beta$, respectively.

\begin{theorem+}{Controlled Resolution Sequences}
Let $P$ and $Q$ be maximal curve systems in $S$. 
There exists a geodesic
$\beta_0,\ldots,\beta_m$ in $\CC_1(S)$ and an elementary move sequence
$P_0\to\ldots\to P_n$, with the following properties:
\begin{enumerate}
\item $\beta_0\in P_0 = P$ and $\beta_m \in P_n = Q$.
\item Each $P_i$ contains some $\beta_j$.
\item $J_\beta$, if nonempty, is always an interval, and if
$[i,j]\subset[0,m]$ then
$$
|J_{\beta_i,\ldots,\beta_j}| \le K (j-i) \sup_Y d_Y(P,Q)^a,
$$
where the supremum is over only those subsurfaces $Y$ whose boundary
curves are components of some $P_k$ with $k\in
J_{\beta_i,\ldots,\beta_j}$. 
\item
If $\beta$ is a curve with non-empty $J_\beta$, then its predecessor
and successor curves $\alpha $ and $\alpha'$ satisfy
$$
| d_\beta(\alpha,\alpha') - d_\beta(P,Q) | \le \delta.
$$
\end{enumerate}
The constants $K,a,\delta$ depend only on the genus of $S$.
The expression $|J|$ for an interval $J$ denotes its diameter.
\end{theorem+}
The sequence $\{P_i\}$ in this theorem is called a resolution
sequence. Such sequences are constructed in
\cite{masur-minsky:complex2} by an inductive procedure: beginning with
a geodesic $\{\beta_i\}$ in $\CC_1(S)$ joining $P$ to $Q$
(we are describing a geodesic here as a sequence of vertices where 
successive ones are joined by edges),
we note that the link of each $\beta_i$ is itself a
curve complex for a subsurface. In each such complex we add a new
geodesic, and repeat. The final structure can then be
``resolved'' into a sequence of maximal curve systems. Control of the
size of the construction at each stage is achieved by applying the
hyperbolicity theorem of \cite{masur-minsky:complex1}.

\subsection{Contraction and quasi-convexity}

Let $\CC(S,\rho,L)$ denote the subcomplex of $\CC$ spanned by the
vertices with $\rho$-length at most $L$. 
We will define a map 
$\Pi_\rho:\CC(S) \to \PP(\CC(S,\rho,L_1))$,
where $\PP(X)$ is the power set of $X$, as follows.  
For $x\in\CC(S)$, let $P_x$ be the curve system associated to the
smallest simplex containing $x$. We define
$$
\Pi_\rho(x) = \bigcup_{f\in\pleat(P_x,\rho)} \short(\sigma_f).
$$
This map turns out to have coarsely the properties of a 
closest-point projection to a convex subset of a hyperbolic space. 

\begin{lemma+}{Contraction Properties}
There are constants $b,c>0$, depending only on the genus of $S$, such
that for any $\rho$ the map 
$\Pi_\rho$ has the following properties: 
\begin{enumerate}
\item (Coarse Lipschitz) If $d_\CC(x,y) \le 1$ then 
$$\diam_\CC (\Pi_\rho(x)\union\Pi_\rho(y)) \le b.$$
\item (Coarse idempotence)  If $x\in\CC(S,\rho,L_1)$ then
$$
d_\CC(x,\Pi_\rho(x)) \le b.
$$
\item (Contraction) If $r=d_\CC(x,\Pi_\rho(x))$ then 
$$
\diam_\CC \Pi_\rho(B(x,cr)) \le b.
$$
\end{enumerate}
\end{lemma+}
Here $d_\CC$ and $\diam_\CC$ refer to distance and diameter measured
in $\CC(S)$, and $B(x,s)$ is a ball of $d_\CC$-radius $s$ around $x$.
By $\Pi_\rho(X)$ for a set $X$ we mean $\union_{x\in X} \Pi_\rho(x)$.

Compare this with the contraction property in
\cite{masur-minsky:complex1}, which was used to prove hyperbolicity of
$\CC(S)$, and the property in \cite{minsky:projections}, which was
used to prove stability properties for certain geodesics in
Teichm\"uller space.

An easy consequence of this theorem is the
following quasiconvexity property for $\CC(S,\rho,L_1)$:
\begin{lemma+}{Quasiconvexity}
If $\beta_0,\ldots,\beta_m$ is a geodesic in $\CC_1(S)$ and
$\beta_0,\beta_m\in \CC(S,\rho,L_1)$, then 
$$
d_\CC(\beta_i,\Pi_\rho(\beta_i)) \le C
$$
for all $i\in[0,m]$ and a constant $C$ depending only on the genus of $S$. 
\end{lemma+}
In particular a geodesic with endpoints in $\CC(S,\rho,L_1)$ never strays too
far from $\CC(S,\rho,L_1)$. This can be compared to the ``Connectivity''
lemma in \cite{minsky:torus}.

The argument for this lemma is very simple, and has its origins in the
stability of quasi-geodesics argument in the proof of Mostow's
Rigidity Theorem \cite{mostow:hyperbolic}: We compare the path
\marginpar{CHECK if Mostow used this...}
$\{\beta_i\}$ to its image ``quasi-path'' $\{\Pi_\rho(\beta_i)\}$. If
the distance between these grows too much then the images slow
down because of the Contraction property (3). Since $\{\beta_i\}$ is a
shortest path and the two paths have nearly the same endpoints
(Coarse idempotence (2)), there is a bound on how far apart they can get.

\medskip

The proof of Lemma \ref{Contraction Properties} is another
application of Thurston's
Uniform Injectivity theorem, as well as some of the tools developed in
\cite{masur-minsky:complex1}. For example, to prove part (1), 
we note that if two vertices of $\CC(S)$ are at distance 1 then they
correspond to disjoint curves, and hence a pleated surface exists that 
maps both geodesically. Thus the argument reduces to bounding
$\Pi_\rho(x)$ for one $x$. Suppose two pleated surfaces share a curve $x$.
If $x$ is short then their short curve sets 
intersect, and we finish by noting that $\diam_\CC(\short(\sigma))$ is
uniformly 
bounded for any $\sigma$. If the curve $x$ is long, then in one of the
pleated surfaces we can find a curve $x'$ of bounded length that runs
along $x$ and then makes a very small jump in its complement
(a long curve in a hyperbolic surface must run very close to
itself). The Uniform  Injectivity theorem is then applied to show that
$x'$ can be realized with bounded length on the second pleated surface
as well.

Part (3) is the main point of the  lemma. Its proof depends
on the analysis in \cite{masur-minsky:complex1}, which 
shows roughly that if $x\in\CC(S)$ is far in $\CC(S)$ from the short
curves of a given hyperbolic metric $\sigma$ on $S$, then sets of the
form $B(x,R)$ for large $R$ 
can be carried in a long nested chain of ``train tracks'' (see
\cite{penner-harer}) 
whose branches mostly run nearly parallel to $x$.
These train tracks are then used to
control $\Pi_\rho(B(x,R))$, via a Uniform Injectivity argument similar
to the previous paragraph.

\subsection{Building a resolution sequence for $\rho$}
\label{finishing argument}

We can now use Theorem \lref{Controlled Resolution Sequences} and
Lemma \lref{Quasiconvexity} 
to produce a resolution sequence adapted to the geometry of our
representation $\rho$. 

As a starting point we need an initial and terminal curve system:
\begin{lemma}{Pants near ends}
  Given $\rho$ with no externally short curves, and a Margulis tube
  $\T_\gamma$ in $N_\rho$, there exist maximal curve systems $P_+$ and
  $P_-$, and pleated surfaces $f_+, f_-$ (in the homotopy class of
  $\rho$) with the following properties:
\begin{enumerate}
\item $P_\pm\in \short(\sigma_{f_\pm})$,
\item $f_+$ and $f_-$ homologically encase $\T_\gamma$.
\end{enumerate}
\end{lemma}

This is done roughly as follows. If $\nu_+$ is a lamination then there
exists a sequence $g_i$ of pleated surfaces exiting the end of
$N_\rho$ corresponding to $\nu_+$. The curves in
$\short(\sigma_{g_i})$ 
converge to $\nu_+$ in the space of laminations (modulo measure), 
and for large enough $i$, $g_i$ can be deformed to $e_+$
without meeting $\T_\gamma$. We can pick $f_+ = g_i$ and let
$P_+\in\short(\sigma_{g_i})$.
The same goes for
$f_-,P_-$, so if both invariants are laminations we have the
conclusion that $f_+$ and $f_-$ must encase $\T_\gamma$.

If the end $e_+$ is geometrically finite we can let $f_+$ be the
pleated map to the convex hull boundary itself, and similarly for
$e_-$. Again choose $P_\pm\in\short(\sigma_{f_\pm}) = \short(\nu'_\pm)$.

\medskip

Note that, if the pleated surfaces $g_i$ are chosen far enough out the
end (in the geometrically infinite case) then the homotopy from $g_i$
to a map in $\pleat(P_+,\rho)$ does not pass through $\T_\gamma$, and
so we may assume $f_\pm\in\pleat(P_\pm,\rho)$ and still have the
encasing condition. When there are geometrically finite ends this is
trickier because $\T_\gamma$ may be close to the convex hull
boundary. Slightly more care is needed in the rest of the construction
in that case. Let us from now on assume that
$f_\pm\in\pleat(P_\pm,\rho)$, and the encasing condition holds.

\medskip

Join $P_+$ to $P_-$ with a resolution sequence 
$P_-=P_0\to\cdots\to P_n=P_+$, as in 
Theorem \ref{Controlled Resolution Sequences}.  Let
$\{\beta_i\}_{i=0}^m$ be the associated geodesic.  This sequence may
be much longer than we need, so we will use Lemma \ref{Quasiconvexity}
to find a suitable subsequence. Recall that we would like our sequence
to have the property of homologically encasing $\T_\gamma$, so let us
try to throw away those surfaces that we are sure cannot meet
$\T_\gamma$.  In particular, let $f\in\pleat(P_i,\rho)$ for some
$i\in[0,n]$, and let $P_i$ contain a curve $\beta_j$.
If $f(S)\intersect \T_\gamma \ne \emptyset$, then
$\gamma$ itself is short in $\sigma_f$  (as in \S\ref{margulis tubes})
and so $\gamma$ is distance 1 from $\Pi_\rho(\beta_j)$. It follows
from Lemma \ref{Quasiconvexity} that 
$$
d_\CC(\beta_j,\gamma) \le C \eqno{(*)}
$$
where $C$ is a new constant depending only on the genus of $S$.
Thus we conclude that
there is a subinterval $I_\gamma$ of $[0,m]$ of diameter 
at most $2C$, such that  $f$ can only meet $\T_\gamma$ when $\beta_j$
satisfies $j\in I_\gamma$. Let us therefore restrict our elementary
move sequence to 
$$P_{s-1}\to\cdots\to P_{t+1}$$
where $[s,t] = \union_{j\in I_\gamma} J_{\beta_j}$, and
renumber it as $P_0\to\cdots\to P_M$.
This subsequence must also encase
$\T_\gamma$, since none of the pieces we have thrown away can meet
$\T_\gamma$. 
Part (3) of  Theorem \ref{Controlled Resolution Sequences} tells us
that 
$$ 
        M \le K(2C)\sup_Y d_Y(P_+,P_-)^a,
$$
where the supremum is over subsurfaces $Y$ whose boundaries appear
among the $P_i$ in our subsequence. This means by ($*$) that the
$\CC(S)$-distance $d_\CC(\boundary Y,\gamma)$ is bounded by $C+1$ for
all such $Y$. The analysis of \cite{masur-minsky:complex2} shows that,
for a fixed such bound, 
$$d_Y(P_+,P_-) \le d_Y(\nu_+,\nu_-) + \delta$$ 
with $\delta $ depending only on the genus of $S$,
provided, when $e_+$ or $e_-$ are geometrically infinite, that the
surface $f_\pm$ are taken sufficiently far out in the ends (for
geometrically finite ends this is an easier consequence of Sullivan's
theorem comparing $\nu_\pm$ with $\nu'_\pm$, though here we must take
a bit more care with the constants to make sure that $\boundary Y$
intersects $P_\pm$).
Since the right side is a priori bounded by hypothesis, we obtain
our desired uniform bound on $M$.

Now let 
$H : S\times[i,i+1]\to N_\rho$ 
be the homotopy provided by Lemma \lref{Elementary Homotopy}, where
$H_i \in \pleat (P_i,\rho)$. After possibly adjusting by
homeomorphisms of $S$ homotopic to the identity, we can piece these
together to a map
$H:S\times[0,M]\to N_\rho$.

Assume first that $\gamma$ is not a component of any $P_i$. 
Then according to Lemma \ref{Elementary Homotopy}, $H$ can make only
uniformly bounded progress through the Margulis tube $\T_\gamma$. Thus
$\diam \T_\gamma$ is bounded above, and $\ell_\rho(\gamma)$ is bounded
below, and we are done.

Now suppose that $\gamma$ does appear in the $\{P_i\}$. Then
$J_\gamma$ is some subinterval of $[0,M]$ by
Theorem \ref{Controlled Resolution Sequences}, and we let $\alpha$ and
$\alpha'$ be the predecessor and successor curves to $\gamma$ in the
sequence. Both of them cross $\gamma$, and we have
by part (4) of Theorem \ref{Controlled Resolution Sequences} that
$d_\gamma(\alpha,\alpha')$ is uniformly approximated by
$d_\gamma(P_+,P_-)$ and hence uniformly bounded.

For simplicity, let us consider now the case that both $\ell_\rho(\alpha)$ and
$\ell_\rho(\alpha')$ are uniformly bounded above and below.
(There is in fact a uniform upper bound on their lengths; 
if they become too short a small additional argument is needed).

Let $\sigma_i \equiv \sigma_{H_i}$ and note that, by Lemma \ref{Elementary
  Homotopy}, for all $i\in J_\gamma$ 
the annuli $\collar(\gamma,\sigma_i)$ coincide. Name this common
annulus $B$. Write $J_\gamma = [k,l]$, and 
consider in $S\times[0,M]$ the solid torus 
$$
U = B\times[k-1,l+1].
$$
By Lemma \ref{Elementary Homotopy}, 
this is the only part of $S\times[0,M]$ that $H$ can map more than a
bounded distance into $\T_\gamma$. The height of this torus, $|k-l+2|$, is
at most $M$ and this is uniformly bounded. The top and bottom annuli
$B\times\{k-1\}$ and $B\times\{l+1\}$ have uniformly bounded geometry
(in $\sigma_{k-1}$ and $\sigma_{l+1}$, respectively),
by the length bounds we've assumed on $\alpha$ and $\alpha'$. 
We will control the size of the {\em meridian} of $U$, and this will
in turn bound the size of $\T_\gamma$.

Assume $\alpha$ is a geodesic in $\sigma_{k-1}$ (where we note its
length is bounded above),
and let $a=\alpha\intersect B$. Similarly
assume $\alpha'$ is a geodesic in $\sigma_{l+1}$ 
and let $a'=\alpha'\intersect B$. The arc $a$ may a priori be
long in $\sigma_{l+1}$,
but its length is estimated by the number
of times it twists around $a'$, or $d_{\AAA(B)}(a,a')$.

A lemma in 2 dimensional hyperbolic geometry establishes
$$
|d_{\AAA(B)}(a,a') - d_\gamma(\alpha,\alpha')| \le C
$$
where this $C$ depends only on $M$, which we have already bounded
uniformly.
The idea of this is that, in each elementary move, the metric $\sigma_i$
changes in a bilipchitz way outside the collars
of the curves involved in the elementary move. From this it follows
that, starting with a geodesic passing through a collar, we obtain a
curve which does only a bounded amount of additional twisting outside
the collar. After $M$ such moves the relative twisting of $\alpha$ and
$\alpha'$ can still be estimated by their twisting inside the collar,
up to an additive bound proportional to $M$.

With this estimate and the bound on $d_\gamma(\alpha,\alpha')$ in
terms of $d_\gamma(P_+,P_-)$, 
we find that $a$
and $a'$ intersect a bounded number of times, so that the length of
$a$ is uniformly bounded in $S\times\{l+1\}$. It follows that the
meridian of $U$
$$
m = \boundary (a\times[k-1,l+1])
$$
has uniformly bounded length in the induced metric. Thus its image is
bounded in $N_\rho$. It therefore spans a disk of bounded diameter,
and in fact we can homotope $H$ on all of $U$ to a new map of bounded
diameter. This bounds the diameter of $\T_\gamma$ from above, and again we
are done.

\section{Conjectures}
\label{conjectures}

\subsection{Length estimates}

The reader may have noticed that in fact the argument
outlined in the previous section shows that the infimum
$\ep=\inf_\gamma\ell_\rho(\gamma)$ and the supremum
$D=\sup_Yd_Y(\nu_+,\nu_-)$ can be bounded one in terms of the other. 
That is,  any positive lower bound for $\ep$ implies some upper bound for
$D$ independent of $\rho$, and vice versa. Thus there is a version of
  the theorem which 
yields non-empty information for quasi-Fuchsian groups (where $\ep>0$
  and $D<\infty$ automatically) as well. However it would be nice to
  have bounds that are more specific and more explicit.

``More specific'' means that we would like to know an estimate on
$\ell_\rho(\gamma)$ for {\em a particular $\gamma$}. In
\cite{minsky:kgcc} we actually show that for any subsurface $Y$, 
a large lower bound on $d_Y(\nu_+,\nu_-)$ implies a small upper bound
for $\ell_\rho(\boundary Y)$. In the other direction something more
complicated would need to be stated, since any curve $\gamma$ can be a
boundary curve for many different subsurfaces. 

``More explicit'' means we would like to know the estimate itself more
explicitly. Furthermore it would be nice to estimate the {\em complex
  translation length} $\lambda$ and not just its real part $\ell$.
In \cite{minsky:torus} this was done for the punctured-torus
case. Here is a possible generalization, stated again in the
case of $\rho$ with no externally short curves. 

\begin{conjecture}{Complex translation length}
Let $\rho$ be a Kleinian surface group with no externally short curves.
There exist $K,\ep>0$ depending only on the genus of $S$ 
such that
$$
\ell_\rho(\gamma)>\ep \implies \sup_{\gamma\subset Y} d_Y(\nu_+,\nu_-)
< K.
$$
Conversely, if $\sup_Y d_Y(\nu_+,\nu_-) \ge K$ then 
$$
\frac{2\pi i}{\lambda_\rho(\gamma)} \asymp
d_\gamma(\nu_+,\nu_-)
+
i 
{\sum^\sim_{\substack{Y\subset S\\
\gamma\subset\boundary Y\\
Y\not\sim\gamma}}} d_Y(\nu_+,\nu_-)
$$
\end{conjecture}
Let us explain the notation used here. 
The expression $\displaystyle\sum^\sim_{x\in X} f(x)$ denotes
$$
1 + \sum_{\substack{x\in X\\ f(x) \ge K}} f(x)
$$
where $K$ is our a-priori ``threshold'' constant.
Our sum then is over all subsurfaces whose boundary contains the
isotopy class of $\gamma$, except for the annulus homotopic to
$\gamma$, excluding those where $d_Y(\nu_+,\nu_-)$ is below $K$.
Both sides of the ``$\asymp$'' symbol are points in the
upper half plane of $\C$, and we take ``$\asymp$'' to mean that the
hyperbolic distance between them is bounded by an a-priori constant
$D_0$. Implicit in the statement is that it holds for some $D_0$
which depends  only on the genus of $S$.

The significance of the hyperbolic distance estimate on 
$2\pi i/\lambda(\gamma)$ is that we can interpret $2\pi
i/\lambda(\gamma)$ as a 
Teichm\"uller parameter for the Margulis tube $\T_\gamma$, as follows
(cf. \cite{minsky:torus} and  McMullen \cite{mcmullen:cusps}).
Normalize $\rho(\gamma)$ so that it acts on  $\hat\C$ by $z\mapsto e^\lambda
z$. The quotient $(\C\setminus\{0\})/\rho(\gamma)$ is then a torus, and there
is a preferred marking of this torus by the pair $(\hat\gamma,\mu)$,
where $\mu$ is the meridian of the torus, or the image of the unit
circle in $\C$, and $\hat\gamma$ is the image of the curve
$\{e^{t\lambda}:t\in[0,1]\}$. Note, this curve depends on the choice
of $\lambda$ mod $2\pi i$. In \cite{minsky:torus} we point out that 
if $\ell_\rho(\gamma)$ is sufficiently short then we can 
choose $\hat\gamma$ to be a minimal representative of $\gamma$ on
the torus just by choosing $\theta = \Im\lambda \in [0,2\pi)$.

The quantity $2\pi i/\lambda$ turns out to be the point in the upper
half-plane representation of the Teichm\"uller space of the torus
which represents the marked quotient torus. Estimating this quantity
up to bounded hyperbolic distance is then equivalent to estimating the
torus structure up to bounded Teichm\"uller distance, which
corresponds to knowing the action of $\rho(\gamma)$ up to bounded
quasi-conformal conjugacy. This finally is equivalent to bounded
bilipchitz conjugacy of the action on $\Hyp^3$, and thus is the
``right'' kind of estimate if we are interested in knowing the
quotient geometry up to bilipschitz equivalence.

The imaginary part of the conjectural estimate is supposed to estimate
the ``height'' of the margulis tube boundary for $\gamma$, and its
real part is supposed to measure the ``twist'' of the meridian around
$\hat \gamma$. In our discussion of the Bounded Geometry Theorem,
we essentially showed that the height was bounded by the
number of elementary moves it took to pass $\T_\gamma$, and the
twisting was bounded by the relative twisting of the predecessor and
successor curves $\alpha$ and $\alpha'$. In general we expect that
large values of $d_Y(\nu_+,\nu_-)$ with $\gamma\subset\boundary Y$
will contribute to parts of the elementary move sequence that make
progress along the sides of $\T_\gamma$, and thus give a good estimate
for its height.

In \cite{minsky:torus} we obtained a similar estimate for the case
where $S$ is a once-punctured torus. (In this case we are not
requiring $S$ to be closed, and our representations must satisfy the
added condition that the conjugacy class corresponding to loops around
the puncture is mapped to parabolics.) Let us state this just in the
case that $\nu_\pm$ are both laminations. For the torus, a 
lamination are determined  by its
slope in $H_1(S,\R) = \R^2$, which takes values in 
$\hhat\R = \R\union\{\infty\}$.
Simple closed curves correspond to rational points.
For any simple closed curve $\alpha$ we defined a quantity
analogous to $d_\alpha(\nu_-,\nu_+)$ as follows: after an appropriate
basis change for $S$ (or equivalently action by an element of
$\SL 2(\Z)$, we may assume that $\alpha$ is represented by $\infty$,
and let $\nu_-(\alpha),\nu_+(\alpha)$ be the irrational numbers
representing the ending laminations. Then define
$$w(\alpha) = \nu_+(\alpha)-\nu_-(\alpha).$$ 
We showed that $\ell_\rho(\alpha)$ can only be short if $w(\alpha)$ is
above a uniform threshold, and  in this case we estimated
$$
\frac{2\pi i}{\lambda_\rho(\alpha)} \asymp w(\alpha) + i.
$$
In fact $w(\alpha)$ is just a measure of relative twisting of $\nu_-$
and $\nu_+$ around $\alpha$, and it is not hard to see that $|w(\alpha)|$
is estimated by our $d_\alpha(\nu_-,\nu_+)$, up to a uniform additive
error. Thus, this is really the same estimate as in Conjecture
\ref{Complex translation length}, since there are no essential subsurfaces
in $S$ other than annuli.

\subsection{General representations}

All the methods that we have presented here depend heavily on the
assumption that $\rho$ is both faithful and discrete. It can be
argued, however, that a full understanding of the deformation space of
hyperbolic structures on a manifold would require some better geometric
description of the whole representation variety, including indiscrete or
non-faithful points, and it is tempting to try to enlist the complex
of curves for this purpose.

The only results I know that offer any hope are in a paper of 
Bowditch \cite{bowditch:markoff}, in which he studies general
representations for the once-punctured torus (again with the
parabolicity condition for the puncture). Such a representation
determines a trace (closely related to complex translation length) for
every conjugacy class, and  in particular for the  simple closed
curves, which in this case correspond to $\Q\union\{\infty\}$, viewed as
the vertices of the Farey tesselation of the disk. To every triangle and
adjacent pair of triangles is associated a relation among the traces
of the vertices, coming from the standard trace identities in $\SL
2(\C)$. Bowditch uses these relations alone, without discreteness, to
analyze the global properties of the trace function, in particular obtaining
a connectedness property for sublevel sets closely analogous to the
quasi-convexity 
property of Lemma \ref{Quasiconvexity}. Using this he is able to
define an invariant of the representation that generalizes the ending
lamination for discrete representations; but it is hard to know how to
extract more information from this invariant.

In the higher genus case, no such analysis has been done, and it would
be very interesting to try it. Elementary moves between pants
decompositions still give rise to trace identities among the curves
involved, although they are a bit more complicated. One wonders at
least whether a result like Lemma \ref{Quasiconvexity} can be
generalized to all representations.

Bowditch is led to the following question: Consider the quantity
$$
\frac{\ell_\rho(\gamma)}{\ell_{\rho_0}(\gamma)}
$$
where $\rho_0$ is some fixed Fuchsian representation, $\rho$ is a
general representation, and $\gamma$ is a non-trivial element of $\pi_1(S)$. 
The infimum of this ratio is positive for quasi-Fuchsian
representations. For a non-quasi-Fuchsian discrete, faithful
representation, the infimum is 0, and can be achieved by considering
only $\gamma$ with {\em simple} representatives. 
The limit points of minimizing sequences in the space of laminations
give the ending laminations for $\rho$.

If $\rho$ is indiscrete or non-faithful the infimum is
again 0 (indeed $\inf\ell_\rho$ is 0 as well), but the question is, 
is the infimum also 0 for the simple elements. In 
other words:

\begin{question} 
Let $S$ be a closed surface of genus at least 2, and let
let $\rho:\pi_1(S) \to \PSL 2(\C)$ be a representation. If
$$
\inf
\frac{\ell_\rho(\gamma)}{\ell_{\rho_0}(\gamma)}
> 0
$$
where $\gamma$ varies over all simple loops in $S$, must
$\rho$ be quasi-Fuchsian?
\end{question}

This question appears to be difficult, and a positive answer
would be a good starting point in using the complex of curves to
analyze general representations.  To indicate its difficulty, note
that it is closely related to the following:

\begin{question}
If $\rho:\pi_1(S) \to \PSL 2(\C)$ is any representation with
non-trivial kernel, does the kernel contain elements represented by
simple loops? 
\end{question}
A positive answer to this question is at least as hard to prove
as the simple loop conjecture for 
hyperbolic 3-manifolds; see Gabai \cite{gabai:simple} and Hass
\cite{hass:simple}. 

\providecommand{\bysame}{\leavevmode\hbox to3em{\hrulefill}\thinspace}


\end{document}